\documentclass[12pt]{article}
\usepackage{amsfonts,amsbsy,amssymb,amsmath,hyperref}
\usepackage{color}
\newcounter{dummy}

\newtheorem{thm}[equation]{Theorem}
\newtheorem{prp}[equation]{Proposition}
\newtheorem{lm}[equation]{Lemma}

\newtheorem{cor}[equation]{Corollary}
\newtheorem{example}[equation]{Example}
\newtheorem{obs}[equation]{Observation}

\newcounter{ilsc}
\newenvironment{ils}{
\begin{list}{\arabic{ilsc}.}{
\topsep10pt
\itemsep10pt
\parsep2pt
\leftmargin0pt
\rightmargin0pt
\listparindent0pt
\itemindent0pt
\labelsep5pt
\labelwidth-5pt
\usecounter{ilsc}
}}
{
\end{list}
}

\def\[#1{\begin{equation}\label{#1}}
\def\]{\end{equation}}

\def\set#1#2{\left\{#1\ \vline\hspace{110sp}\vline\hspace{110sp}\vline\ #2\right\}}

\def\Z{\mathbb{Z}}

\def\R{\mathbb{R}}

\def\RP{\mathbb{R}\mathrm{P}}
\def\CP{\mathbb{C}\mathrm{P}}
\def\HP{\mathbb{H}\mathrm{P}}
\def\P{\mathbb{P}}
\def\U{\mathrm{U}}
\def\SU{\mathrm{SU}}

\def\Span{\mathrm{Span}}
\def\r){\right)}
\def\l({\left(}

\def\diag#1{\mathrm{diag}\left(#1\right)}

\def\sq{\hbox{\rlap{$\sqcap$}$\sqcup$}}
\def\e{\epsilon}
\def\vk{((1+\e)\nu/f^2)v^{\ast}_k,v^{\alpha}_k)}

\def\la{\langle}
\def\ra{\rangle}
\def\fra{\frac{(1+\epsilon)\nu}{f^2}}
\def\quo{D^{n+1}\times_{\alpha} G/K}
\def\pro{D^{n+1}\times G/K}
\def\r{\mathrm{Ric}}
\def\pf{\noindent {\bf Proof. }}
\def\pfo#1{\medskip\par\noindent {\bf Proof of #1}}
\def\proofend{\hfill\sq\par}
\def\den{(1+\epsilon)\nu+f^2}

\begin{document}\sloppy

\title{On the curvature on $G$-manifolds with finitely many non-principal orbits}
\author{S. Bechtluft-Sachs, D. J. Wraith}
\date{}
\maketitle

\begin{abstract}
{\small We investigate the cuvature of invariant metrics on $G$-manifolds with finitely many non-principal orbits.
We prove existence results for metrics of positive Ricci curvature and non-negative sectional curvature, and discuss some families of examples to which these existence results apply.}
\end{abstract}


\section{Introduction}
In this paper we will study the geometry of $G$-manifolds with finitely many non-principal orbits. Here, both the group $G$ and the manifold are smooth, compact and connected, and the action of $G$ on the manifold is smooth and effective.
The orbits of such a $G$-action are either principal, exceptional (that is, non-principal but with the same dimension as a principal orbit), or singular. The codimension of the principal orbits is the cohomogeneity of the $G$-space.

The motivation for studying the situation where the non-principal orbits are finite in number arises from the study of cohomogeneity-one manifolds where existence of invariant metrics with positive Ricci curvature is controlled by the fundamental group.
\begin{thm}[\cite{Be},\cite{GZ2}] A compact $G$-manifold of cohomogeneity $0$ or $1$ admits an invariant metric of positive Ricci curvature if and only if its fundamental group is finite.
\end{thm}
It is already observed in \cite{GZ2} that this can not carry over to cohomogeneity $\geq 4$. The situation in between, i.e. for cohomogneity $2$ and $3$ is essentially open. There are some partial results however, under stronger conditions. Thus in \cite{BW0} metrics of positive Ricci curvature are constructed on asystatic compact $G$-manifolds with finite fundamental group all of whose singular orbits are fixed points.

Cohomogeneity-one manifolds have been studied intensively in recent times. The reason that these objects form such a good family to study is that they have a simple topoplogical description, but form a large and rich class containing many interesting and important examples. A compact cohomogeneity one manifold is either a fibre bundle over a circle (in which case all orbits are principal), or has precisely two non-principal orbits. Of particular note is the role that cohomogeneity one manifolds continue to play in the search for new examples of manifolds with good curvature characteristics. If one considers invariant metrics, then symmetry reduces the problem of describing and analysing such metrics to one which has a reasonable chance of being tractable. For example, new families of manifolds with non-negative sectional curvature, including many exotic spheres in dimension seven, have been been discovered as a result of this approach \cite{GZ1}. The cohomogeneity one condition in the context of positive sectional curvature has attracted particular attention due to the work of Grove, Ziller, Wilking, Verdiani and others. (See for example \cite{GWZ}, \cite{GVWZ}, \cite{V1}, \cite{V2}.) In a recent development, Grove, Verdiani and Ziller \cite{GVZ} and independently Dearricott \cite{D} have announced the existence of a new cohomogeneity-one manifold with positive sectional curvature. Together with the recent announcement of a positive sectional curvature metric on the Gromoll-Meyer sphere by Petersen and Wilhelm \cite{PW}, these are the first new examples of manifolds admitting positive sectional curvature metrics for a number of years.

As a compact cohomogeneity one manifold has either zero or two non-principal orbits, the study of $G$-manifolds with finitely many non-principal orbits can be viewed as a natural generalisation of the cohomogeneity one situation.

\medskip
Topologically, $G$-manifolds with finitely many non-principal orbits have a rather simple structure. Let $M$ be such a $G$-manifold, and suppose the principal isotropy is $K$, so that principal orbits are all equivariantly diffeomorphic to $G/K$. Suppose the non-principal isotropy groups are $H_1,...H_p$, so the non-principal orbits are $G/H_1,...,G/H_p$. It is crucial that the non-principal isotropy groups $H_i$ act with only one orbit type on spheres. In cohomogeneity one their action is transitive. In higher cohomogeneity, $K$ is normal in $H_i$ and $L_i:=H_i/K$ acts freely on a sphere. By Theorem 6.2 in \cite{Br} the group $L_i$ must finite or one of the groups $\U(1)$, $\SU(2)$, $N_{\SU(2)}\U(1)$. For details on this and the subsequent topological facts about $G$-manifolds with finitely many non-principal orbits we refer to \cite{BW}. From now on we will assume that the cohomogeneity is at least $2$.

Let $N_1,...N_p$ be disjoint equivariant tubular neighbourhoods of the non-principal orbits. Then $M^0:=M-\cup_{i=1}^p N_i$ consists of principal orbits only and we have a $G/K$-fibre bundle $M^0\to B$ with structure group $W:=N_GK/K$. Each tubular neighbourhood $N$ admits a simple description. Let $T:=\partial N$. It is clear that $T$ is a sphere bundle over $G/H$. Let $L$ denote one of the groups $\U(1)$, $\SU(2)$, $N_{\SU(2)}\U(1)$, or a finite group $\Gamma\subset O(n+1)$ which acts freely on $S^n$ and $\alpha\colon L \rightarrow H/K$ be an isomorphism. This naturally defines an action of $L$ on $D^{n+1} \times G/K$ (where $z\in L$ sends $(x,gK)\mapsto (zx,g\alpha(z^{-1}))K$). We will use the symbol $\times_{\alpha}$ to indicate quotients under this action. Thus we have $T\cong S^n\times_{\alpha} G/K$ and $N\cong D^{n+1}\times_{\alpha} G/K$.

To each non-principal orbit corresponds a boundary component of $B$, $\partial_i B=S^{n_i}/L_i.$ Since $\partial_i B$ is a quotient of a sphere by a free action, it follows that in the case of a singular orbit $\partial_i B$ is a quaternionic projective space if $L_i=\SU(2)$, a complex projective space if $L_i=\U(1)$, and $\CP^{odd}/\Z_2$ if $L_i=N_{\SU(2)}\U(1)$. If the orbit is exceptional, we must have $L_i$ finite, and in odd cohomogeneity $\partial_i B=\RP^{n_i}$, $n_i$ even, and $L_i=\Z_2$. Notice that if there is a singular orbit, the cohomogeneity must be odd.

In contrast to the case of cohomogeneity one, any number of non-principal orbits can occur. Thus, for instance, there are actions of $\U(1)$ on $S^{2k+1}$ and of $\SU(2)$ on $S^{2k}$ with only one non-principal orbit. We will show that many of the examples of $G$-manifolds with finitely many non-principal orbits constructed in \cite{BW} admit an invariant metric with positive Ricci curvature. We recall a few examples:
\begin{example}\label{list} (`Doubles') Let $L$ be finite or one of $\U(1)$, $\SU(2)$, $N_{\SU(2)}\U(1)$. Let $G$ and $K$ be compact Lie groups, $K\subset G$ and $\alpha\colon L\to N_GK$ be injective. Let
$$M:= D^{n+1}\times_{\alpha} G/K\cup D^{n+1}\times_{\alpha} G/K$$
where we glue the common boundary $T=S^n\times_{\alpha} G/K$ via the identity. This
is a $G$ manifold with two identical singular orbits. The orbit space is the suspension of $S^n/L$.
\end{example}
Recall that for $p_1,p_2$ coprime, the Aloff-Wallach space $W_{p_1,p_2}$ is the quotient $\SU(3)/\set{\diag{z^{p_1},z^{p_2},z^{-p_1-p_2}}}{z\in\U(1)}$. 
\begin{example}\label{dim11} Given any two Aloff-Wallach spaces $W_{p_1,p_2}$ and $W_{q_1,q_2}$, there is an 11-dimensional
$SU(3)$-manifold $M^{11}_{p_1p_2q_1q_2}$ of cohomogeneity three, orbit space $S^3$, and two singular orbits equal to the given Aloff-Wallach spaces.
Within this family there is an infinite sequence of pairwise non-homotopy equivalent manifolds for which each pair of singular orbits is non-homotopy equivalent. There is also an infinite sequence of pairwise non-homotopy equivalent `doubles', that is, manifolds with two identical singular orbits.
\end{example}
\begin{example}\label{dim13} Given Aloff-Wallach spaces $W_{p_1,p_2}$ and $W_{q_1,q_2}$, there is a 13-dimensional $\SU(3)$-manifold
$M^{13}_{p_1p_2q_1q_2}$ of cohomogeneity 5, orbit space $\Sigma \CP^2$, and two singular orbits equal to the given Aloff-Wallach manifolds if and only if
$p_1^2+p_1p_2+p_2^2=q_1^2+q_1q_2+q_2^2$. Within this family there is an infinite
sequence of pairwise non-homotopic manifolds for which each pair of singular orbits is non-homotopic. There is also
an infinite sequence of pairwise non-homotopic `doubles', that is, manifolds for which each pair of singular orbits is identical.
\end{example}

We now consider the geometry of invariant metrics on $G$-manifolds with finitely many non-principal orbits.
In \cite{LY} it was shown that the existence of a compact non-abelian Lie group action on a compact manifold means that the manifold admits a metric of positive scalar curvature. Moreover, the same construction actually yields an {\it invariant} metric of positive scalar curvature (a fact not pointed out in \cite{LY}, but observed, for example, in \cite{H}).
Thus if $M$ is a $G$-manifold of the type under consideration in this paper, provided $G$ is not a torus, then $M$ must admit an invariant metric of positive scalar curvature. \par

Of particular importance here is the fact that for cohomogeneity one, the space of orbits is one dimensional and so makes no contribution to the curvature. In higher cohomogeneities, this is no longer the case.
Indeed the space of orbits might have particularly bad Ricci curvature characteristics. Thus there seems little hope of being able to prove a positive Ricci curvature existence theorem of comparable generality to \cite{GZ2} in our situation.
It seems reasonable to expect that the Ricci curvature of the space of orbits will play an explicit role in any existence theorem.
In fact, we are able to prove the following:
\begin{thm}\label{thm-riccurvmain} Let $M$ be a compact $G$-manifold with finitely many singular orbits, for which the principal orbit $G/K$ has finite fundamental group.
Let $M^0$ be the manifold with boundary resulting from the removal of small invariant tubular neighbourhoods around the non-principal orbits, so $M^0$ is the total space of a $G/K$-bundle with base $B$. The boundary components of $B$ are all quotients of spheres by free actions of subgroups of the orthogonal group, and thus have a standard metric induced by the round metric of radius one.
If $B$ can be equipped with a Ricci positive metric such that
\begin{enumerate}
\item\label{cond-fbs} for each $i$, the metric on boundary component $\partial_i B$ is the standard metric scaled by a factor $\lambda_i^2$;
\item\label{cond-pricurv} the principal curvatures (with outward normal) at boundary component $\partial_i B$ are strictly greater than $-1/\lambda_i$;
\end{enumerate}
then $M$ admits a $G$-invariant metric with positive Ricci curvature.
\end{thm}
From this it is easy to deduce
\begin{cor}\label{cor-4.2}
All compact G-manifolds with two singular orbits, orbit space a suspension of either a projective space or $\CP^{odd}/\Z_2$ and principal orbit $G/K$ with $\pi_1(G/K)$ finite, admit invariant Ricci positive metrics.
\end{cor}
In particular we have:
\begin{cor}\label{cor-4.3}
The families $M^{11}_{p_1p_2q_1q_2}$ of example \ref{dim11}, and $M^{13}_{p_1p_2q_1q_2}$ of example \ref{dim13} all admit invariant metrics of positive Ricci curvature.
\end{cor}
We also have examples of manifolds with a single non-principal orbit and positive Ricci curvature:
\begin{thm}\label{one} For every $n\ge 2$ there is an $SU(n)$-manifold of dimension $n^2+2$ and cohomogeneity three with a single singular orbit and an invariant metric of positive Ricci curvature.
\end{thm}
We now turn our attention to the sectional curvature, and specifically non-negative sectional curvature.
It is not difficult to show that in certain special circumstances, we can obtain manifolds with invariant metrics of non-negative sectional curvature.
\begin{thm}\label{K} All $G$-manifolds with two identical singular orbits and orbit space a suspension of either a projective space or $\CP^{odd}/\Z_2$ admit invariant metrics with non-negative sectional curvature.
\end{thm}
We immediately obtain:
\begin{cor} There are infintely many homotopy types of manifolds in both the families $M^{11}_{p_1p_2q_1q_2}$ and $M^{13}_{p_1p_2q_1q_2}$ which admit invariant metrics of non-negative sectional curvature.
\end{cor}
This paper is laid out as follows. In section 2 we investigate the geometry of tubular neighbourhoods of non-principal orbits.
In section 3 we give proofs for the main results, with some of the more technical results postponed until section 4. We conclude with a collection of open problems in section \ref{sec-open}.

The authors would like to thank Dmitri Alekseevsky for encouraging us to study manifolds with finitely many singular orbits, and for his subsequent help. We would also like to thank Thomas P\"uttmann for reading a preliminary draft of this paper and for his valuable comments.

\section{Tubular neighbourhoods of non-principal orbits}

With the same notation as in the Introduction, let us first focus on the `regular' part $M^0$ of our $G$-manifold $M$. We will assume that
$B=M^0/G$ comes equipped with a Ricci positive metric satisfying the properties given in Theorem \ref{thm-riccurvmain}.
Now fix a bi-invariant
metric $g_0$ on $G$, and let $\nu>0$. The metric $\nu g_0$ induces a normal homogeneous metric (which we will also denote
$\nu g_0$) on $G/K$. By assumption, $\pi_1 G/K$ is finite, and it is well-known (see \cite{Be}) that the metric $\nu g_0$ on $G/K$
has positive Ricci curvature.

On $M^0$ we will introduce a submersion metric. (See \S9 of \cite{B} for more details about the construction of such metrics.)
For this we require three ingredients: a base metric, a fibre metric, and a horizontal distribution of subspaces.
In the current situation we have base and fibre metrics. Any choice of horizontal distribution then gives a
submersion metric, and it follows from (\cite{B}; 9.70) that this submersion metric will have positive Ricci curvature provided
the constant $\nu$ is chosen sufficiently small. From now on, we will assume that $M^0$ is equipped with such
a Ricci positive submersion metric. We are free, of course, to a smaller value of $\nu$ later on if required.

Let us now turn our attention to the tubular neighbourhoods of non-principal orbits. As discussed in the
Introduction, these all take the form
$$N:=D^{n+1} \times_{\alpha} G/K,$$ where $\alpha:L \rightarrow H/K$ is an isomorphism, $L=\U(1),N_{\SU(2)}\U(1),\SU(2)$ or
a finite subgroup of $O(n+1)$, and $H$ is the singular isotropy.

Our approach to constructing a metric on this neighbourhood is to define a metric $g_1$ on $D^{n+1}$, a metric $g_2$ on $G/K,$
and then consider the product metric $g_1+g_2$ on $\pro$. By making a suitable choice of $g_1$ and $g_2$, we can arrange for this product
metric to induce a well-defined metric $g_Q$ on the quotient $N$. Moreover, by a possibly more refined choice of starting metrics,
we can show that this induced metric can always have positive Ricci curvature. Of course, such
neighbourhoods must then be glued smoothly, and within positive Ricci curvature, into $M^0$ with its submersion metric. In particular,
this means that the $G/K$-fibres on the boundary of the tubular neighbourhood must have normal homogeneous metrics $\nu g_0.$

We first deal with the case of exceptional orbits.

\begin{thm}\label{ex} Consider a tubular neighbourhood $\quo$ of an isolated exceptional orbit $G/H$, where $\alpha:\Gamma \rightarrow H/K$
is an isomorphism from a finite subgroup $\Gamma \subset O(n+1).$
Fix a bi-invariant background metric $g_0$ on $G$, and let $\nu>0$. Given constants $\lambda>0$ and $0<\Lambda<1$,
there is a $G$-invariant Ricci positive metric $g_Q$ on $\quo$ such that the complement of the exceptional orbit has a submersion metric
with fibres isometric to $(G/K,\nu g_0)$ and base isometric to $((0,R] \times S^n/\Gamma, dr^2 + h^2(r)\sigma^2)$,
where $\sigma^2$ is the metric on $S^n/\Gamma$ induced by $ds^2_n$, and where $h(R)=\lambda$ and $h'(R)=\Lambda.$
\end{thm}

\pf Choose a function $h(r)$ such that $h(r)=\sin r$ for $r$ small, $h''(r)<0$ for all $r$, and $h(R)=\lambda$, $h'(R)=\Lambda$ for some $R>0.$
It is clear that we can make such a choice. Moreover, the resulting metric $dr^2+h^2(r)ds^2_n$ on $D^{n+1}$ will have positive Ricci curvature, and hence
so will the product metric $dr^2+h^2(r)ds^2_n+\nu g_0$ on $D^{n+1} \times G/K$.  As $\Gamma$ acts isometrically on this product,
we obtain a well-defined metric $g_Q$ on the quotient $\quo$. As the quotient map is a finite covering, this induced metric is locally
isometric to $dr^2+h^2(r)ds^2_n+\nu g_0$, and hence has positive Ricci curvature. Moreover, the metric on each
$G/K$ fibre in $(D^{n+1}-\{0\}) \times_{\alpha} G/K$ is clearly isometric to $\nu g_0.$
\proofend
\smallskip
We next consider the case where $L=\U(1)$ or $\SU(2).$

\begin{thm}\label{tube} Consider a tubular neighbourhood $\quo$ of an isolated singular orbit $G/H$. Fix a bi-invariant background
metric $g_0$ on $G$ so that $\alpha:L \rightarrow H/K$ is an isometry, where $L=\U(1)$ or $L=\SU(2).$
Given constants $\lambda>0$, $0<\Lambda<1$ and $0<\e<\e_0$ (where $\e_0=\e_0(G,H,K,g_0)$ is
the constant from Observation 5.5), for any $\nu<\lambda\e/(1+\e),$
there is a $G$-invariant Ricci positive metric $g_Q$ on $\quo$ such that for some small $\iota>0$, the $\iota$-neighbourhood
of the boundary has a submersion metric with fibre $G/K$, base $[R-\iota,R] \times \P$, all fibres isometric to the normal
homogeneous metric induced by $\nu g_0$, and base metric $dr^2+h^2(r)ds^2_n$ with $h(R)=\lambda$ and $h'(R)=\Lambda.$
\end{thm}

The proof of this result is somewhat technical, and depends crucially on several explicit curvature calculations. For this reason, we postpone
the proof and the relevant computational lemmata until section 4.

It remains to consider the case when $L=N_{\SU(2)}\U(1).$

\begin{cor}\label{norm} The statement of Theorem \ref{tube} continues to hold in the case $L=N_{\SU(2)}\U(1).$
\end{cor}

\pf The isomorphism $\alpha:N_{\SU(2)}\U(1) \rightarrow H/K$ restricts to an isomorphism $\alpha_0:U(1) \rightarrow (H/K)_0$ where
$(H/K)_0$ is the identity component of $H/K$. By Theorem \ref{tube} we obtain a Ricci positive metric on $D^{n+1} \times_{\alpha_0} G/K$
with the desired properties. We now simply observe that $\quo$ is a $\Z_2$-quotient of $D^{n+1} \times_{\alpha_0} G/K$, with $\Z_2$
acting isometrically. Thus we obtain a Ricci positive metric on $\quo$, and it is clear that this metric has all the claimed properties.
\proofend

\section{Proofs of the main results}\label{pf}

\noindent {\bf Proof of Theorem \ref{thm-riccurvmain}.} We make use of an observation due to Perelman \cite{P}, both to control the form of the metric on the space of orbits $B$ near the boundary components, and then to allow smooth Ricci positive gluing with the tubular neighbourhoods of isolated singular orbits.
According to Perelman, given two Ricci positive manifolds with isometric boundary components, if the principal curvatures at one boundary component are (strictly) greater than the negatives of the corresponding principal curvatures for the other boundary component, the non-smooth metric which results from gluing the two boundary components together can be smoothed in an arbitrarily small neighbourhood of the join to produce a metric with global Ricci positivity. \par
With this in mind, we add small collars of the form $\P_i\times [0,\e_i]$ for some small $\e_i$ to each boundary component $\P_i$.
On this collar we assume a metric of the form $ds^2+\theta_i(s)g_{\P}$ with $g_{\P}$ the standard metric on $\P$.
In order for the boundary $\P_i$ to be isometric to the collar at $s=0$ we clearly need $\theta_i(0)=\lambda_i$.
If we choose $\theta_i$ such that $\theta_i''<0$ and $|\theta_i'|<1$, it is easy to check that the Ricci curvature of the collar will be positive.
Let $p_i$ be the infimum of the principal curvatures (with outward pointing normal) at the boundary component $\P_i$.
Recall that by the assumption (2) in the statement of the Theorem, $p_i>-1/\lambda_i$.
It is easy to see that the principal curvatures at the $s=0$ boundary component of the collar (again with outward normal) are all equal to $-\theta_i'(0)/\theta_i(0)$.
Therefore if $\theta_i'(0)<\lambda_i p_i$ for all $i$, then by Perelman we can join the collar to the space of orbits and smooth (in an arbitrarily small region) within Ricci positivity.
(Note that such a value of $\theta_i'(0)$ always exists: by assumption, $p_i>-1/\lambda_i$, so $p_i\ge -(1/\lambda)+c_i$ for some $c_i>0$.
Therefore the upper bound on $\theta_i'(0)$ is $-1+\lambda_ic_i>-1$, and so we have a well-defined non-empty interval $(-1,-1+\lambda_i c_i)$ from which to choose $\theta_i'(0)$.) Notice that by making a careful choice for $\theta$ and $\e_i$, we can also ensure that at the new boundary we create, the metric still satisfies conditions (1) and (2) in the statement of the Theorem.
Specifically, we need $\theta_i'(\e_i)>-1$ in order to satisfy these requirements.
The conclusion from this analysis is that without loss of generality we are free to make the following
\par\noindent{\it Assumption.} The metric in a small
neighbourhood of the boundary $\P_i$ is isometric with $(\P_i \times [0,\e_i];ds^2+\theta_i(s)g_{\P})$ for some $\e_i>0,$
and with all principal curvatures at the boundary equal to $\theta_i'(\e_i)/\theta_i(\e_i) \in (-1/\theta_i(\epsilon),0).$
\par
Now let us turn our attention to $M$ itself.
Recall that $M^0$ is the manifold with boundary resulting from the removal of small invariant tubular neighbourhoods $N_i$ around the non-principal orbits,
so $M^0$ is the total space of a $G/K$-bundle over a manifold with boundary $B$.
Consider the tubular neighbourhood $\quo$, and suppose we wish to glue this to boundary component $i$ of $M^0$. Let $\P$ be the corresponding
boundary component of $B.$
Equip $\quo$ with the metric $g_Q$ as in Theorem \ref{tube} or Corollary \ref{norm} in the singular cases, or as in Theorem \ref{ex} in the exceptional case, and where we have chosen $\lambda=\theta_i(\e_i)$ and $\Lambda=|\theta_i'(\e_i)|$, with $\theta_i$ as in the above assumption.
The scaling function $h$ in $g_Q$ will now glue with $\theta$ when the $r$ and $s$ parameters are suitably concatenated to create a $C^1$ scaling function.
As a pre-requisite for the smooth gluing of $\quo$ to $M^0$, we need this scaling function to be smooth.
We can easily achieve this by making a minor adjustment to $h(r)$ close to the boundary of $\quo$, and in particular for $r$ in the interval $(R-\iota,R]$ (with $\iota$ as in Theorem \ref{tube}) in the singular case.
Specifically, we can adjust $h$ to make the required alteration in the second derivative, whilst keeping the variation in both $h$ and $h'$ arbitrarily small.
Such an adjustment does not destroy positive Ricci curvature: this is easy to see, for instance by using Proposition \ref{prp-5.4}. \par
Now consider the boundary of $\quo$ and the corresponding boundary of $M^0$ as $G/K$-bundles over $\P$.
A further pre-requisite for smooth metric gluing is that the horizontal distributions of these bundles much match under the identification.
We are free to choose a horizontal distribution for $M^0$ (viewed as as $G/K$-bundle over $B$), so we choose this in a way so as to agree near the appropriate boundary component with that coming from $g_Q$.
Note that the horizontal distribution determined by $g_Q$ is induced from the standard horizontal distribution for the Hopf fibration
in the singular case and from $TS^n$ in the exceptional case. Thus it is independent of all choices involved in the construction of $g_Q$, and
in particular is independent of the functions $f$, $h$ and the constant $\nu$.
This last point is important, as it means we can select a horizontal distribution for $M^0$ (taking all boundary components into consideration) at the outset, and thus find a constant $\nu_0>0$ such that the submersion metric on $M^0$ has positive Ricci curvature for all $\nu<\nu_0$. \par
It remains consider the $G/K$-fibres on both sides of the join.
The metric on the fibres near the boundary of $\quo$ is, by the construction of $g_Q$, the normal homogeneous metric induced by $\nu g_0$ on $G$, and this agrees with the fibre metrics for $M^0$.
According to Theorem \ref{tube}, in the singular case we need $\nu<\lambda_i\e/(1+\e)$ for the Ricci positivity of tubular neighbourhood $N_i$, and for Ricci positivity of $M^0$ we need $\nu<\nu_0$ as above.
Therefore, our construction of a smooth global Ricci positive metric can be completed by choosing a value for $\nu$ which is less than $\nu_0$
in the exceptional case, and less than both $\nu_0$ and the minimum of the $\lambda_i\e/(1+\e)$ in the singular case.
\proofend

\pfo {Theorem \ref{one}.}
We form a tubular neighbourhood of the singular orbit by setting $N=D^4\times_{\alpha} SU(n)$ where $\alpha:U(1)\rightarrow SU(n)$ is any injective homomorphism.
Since $SU(n)$ is simply-connected, the boundary of this neighbourhood $\partial N=S^2\times SU(n)$.
We can therefore equivariantly glue a product $D^3\times SU(n)$ to this neighbourhood to produce a closed manifold.
The existence of an invariant Ricci positive metric now follows easily from Theorem \ref{thm-riccurvmain}.
The manifold $B$ in the statement of Theorem \ref{thm-riccurvmain} is simply $D^3$ in our situation, and thus the metric conditions on $B$ required by Theorem \ref{thm-riccurvmain} can easily be satisfied.
\proofend

\pfo {Theorem \ref{K}.}
Construct a metric $g$ on the tubular neighbourhood $N=\quo$ of a singular orbit as the quotient of a product metric $g_1+g_2$ on $\pro$ in the following way.
Suppose $D^{n+1}$ has radius $\pi/2$ and set $g_1=dr^2+\sin^2 r\ ds^2_n$.
(So $g_1$ is round with constant sectional curvature $1$.) Let $g_2$ be a normal homogeneous metric on $G/K$.
The product metric $g_1+g_2$ clearly has non-negative sectional curvature as the curvatures of both $g_1$ and $g_2$ are non-negative.
As Riemannian submersions are non-decreasing for the sectional curvature, it follows that $g$ also has non-negative sectional curvature.
The $G/K$ fibre metrics at the boundary of $N$ are not normal homogeneous, but this does not matter in the case of doubles as we wish to glue $N$ to an identical object.
Viewing a neighbourhood of $\partial N$ as an $S^n$-bundle over $G/H\times (\pi/2-\epsilon,\pi/2]$, the only issue we need to consider when gluing is the smoothness of the $S^n$ metrics across the join.
But this is clear as $\sin r$ for $r\in [0,\pi/2]$ concatenates smoothly with its `reverse' $\sin ((\pi/2)-s)$ for $s\in [0,\pi/2]$ when $r=\pi/2$ is identified with $s=0$.
\proofend

\section{Curvature computations}\label{sec-proof} In section 3 we gave the proof of Theorem \ref{thm-riccurvmain}, our main existence result for positive Ricci curvature. This proof depends crucially on Theorems \ref{ex}, \ref{tube} and Corollary \ref{norm}. However, as a result of its technical nature, the proof of
Theorem \ref{tube} was postponed. The aim of the current section is to establish this Theorem and to perform the pre-requisite curvature computations.
In order to do this, we must study metrics on
$$D^{n+1}\times_{\alpha} G/K,$$
where in this case $\alpha$ is a group isomorphism $L \rightarrow H/K$ with $L=\U(1)$ or $\SU(2)$. In the sequel, it will be convenient
to identify $\U(1)$ and $\SU(2)$ with the spheres $S^1$ respectively $S^3$. We will generally deal with both cases at once by writing $S^q.$
Recall that the group $S^q$ acts on the product $\pro$ as follows:
$$z(p,gK)=(zp,gK\alpha(z^{-1})),$$
where the $S^q$-action on $D^{n+1}$ is the standard Hopf action on the first factor of $(S^n\times [0,R])/(S^n\times\{0\})=D^{n+1}$, and
the expression $\quo$ denotes the quotient of $D^{n+1}\times G/K$ by this action. \par
Let $\{v_k\}$ denote a local orthonormal frame field for $S^q(1)$.
As $S^q$ acts via $\alpha$ on $G/K$ and directly on $S^n$, we obtain induced action fields which we will denote $\{v^{\alpha}\}$ and $\{v^{\ast}\}$ respectively.
Notice that $\{v^{\alpha}\}$ is a local frame field for $H/K$-orbits in $G/K$, and that $\{v^{\ast}\}$ is a local frame field for the fibres of the Hopf fibration. \par
We will construct a product metric $g_1+g_2$ on $\pro$ in such a way that it induces a well-defined Ricci positive $G$-invariant metric on $\quo$, with all the properties we will need to glue smoothly into $M^0$.
First of all, we concentrate on constructing $g_1$ for the $D^{n+1}$-factor.
It is technically easier if we remove the centre point of $D^{n+1}$ and view the space as $(0,R]\times S^n$ for some $R>0$.
We will construct metrics on $(0,R]\times S^n$, but the boundary conditions we impose will ensure that our metric extends smoothly to $D^{n+1}$.
\par
Let $g_{\P}$ denote the standard Fubini-Study metric on a complex or quaternionic projective space $\P$.
In the Lemma below, we make use of the O'Neill formulas for the Ricci curvature of a Riemannian submersion.
See \cite{B} \S9 for details of the formulas and definitions of the terms involved.
See \cite{B} \S9.59, or \cite{V} for a discussion about constructing submersion metrics.
\begin{lm}\label{lm-5.1}
Consider the `extended' Hopf fibration $S^q\hookrightarrow S^n\times (0,R]\rightarrow\P\times (0,R]$.
Equip the base with the metric $dr^2+h^2(r)g_{\P}$ and the fibres with the metric $f^2(r)ds^2_q$.
Introduce into the total space the horizontal distribution which is the obvious extension of the stadnard horizontal distribution for the Hopf fibration.
Let $g_1=g_1(f,h)$ denote the resulting submersion metric.
Then $g_1$ has the following Ricci curvatures (denoted $\r_1$):
\begin{align*}
&\r_1(\partial_r)=-\dim\P\frac{h''}{h}-q\frac{f''}{f}\\
&\r_1(X_i)=\frac{1}{h^2}\r_{\P}(\check{Y}_i)-\frac{h''}{h}-q\frac{f'h'}{fh}-2\frac{f^2}{h^4}\la A_{Y_i},A_{Y_i}\ra\\
&\r_1(v^{\ast}_k)=\frac{1}{{f^2}}\r_{S^q(1)}(v_k)-\frac{f''}{f}-(q-1)\left(\frac{f'}{f}\right)^2-\dim\P\frac{f'h'}{fh} +\frac{f^2}{h^4}\la Av^{\ast}_k,Av^{\ast}_k\ra\\
&\r_1(X_i,v^{\ast}_k)=\frac{f}{h^3}\la (\check{\delta} A)Y_i,v^{\ast}_k\ra
\end{align*}
where $\{X_i\}$ are an orthonormal spanning set of vector fields for the horizontal distribution on $S^n$, $Y_i=h^{-1}X_i$ and $\check{Y}_i$ is the projection of $Y_i$ to $\P$.
The terms $\la A_{Y_i},A_{Y_i}\ra$, $\la Av^{\ast}_k,Av^{\ast}_k\ra$ and $\la (\check{\delta} A)Y_i,v^{\ast}_k\ra$ are the standard terms for the Hopf fibration on $S^n(1)$.
All other mixed Ricci curvature terms vanish.
\end{lm}
\pf These formulas follow from \cite{W2}, Proposition 4.2.
\proofend
We can make the above formulas more explicit by substituting the appropriate values for $\r_{\P}$, $\r_{S^q(1)}$ and the $A$-tensor terms.
It is well known \cite{B}, p.258, that with their standard Fubini-Study metrics, both $\CP^m$ and $\HP^m$ are Einstein manifolds, with Einstein constants $2m+2$ and $4m+8$ respectively.
For the $A$-tensor terms we have:
\begin{lm}\label{lm-5.2}
For the standard Hopf fibration over $\CP^m$ or $\HP^m$ we have $\la Av^{\ast}_k,Av^{\ast}_k\ra$ equal to $\dim \CP^m$ respectively $\dim \HP^m$, $\la AY_i,AY_i\ra=q$ and $\la (\check{\delta} A)Y_i,v^{\ast}_k\ra=0$.
\end{lm}
\pf These expressions can be evaluated by applying the O'Neill formulas for the Ricci curvature \cite{B}, \S9.70, to the standard Hopf fibration, using the known constant Ricci curvature values for the base, total space and fibre.
The computations are all elementary.
We mention only $\la (\check{\delta} A)Y_i,v_k\ra$ in the case where the base is $\CP^m$.
For this, consider $\r(Y_i+v^{\ast}_k)=2(2m)$ as $\|Y_i+v^{\ast}_k\|=\sqrt{2}$.
We also have
\begin{align*}
\r(Y_i+v^{\ast}_k)&=\r(Y_i)+\r(v^{\ast}_k)+2\r(Y_i,v^{\ast}_k)\\
&=2(2m)+2\r(Y_i,v^{\ast}_k).
\end{align*}
Hence $\r(Y_i,v^{\ast}_k)=0$.
The O'Neill formulas show that $\r(Y_i,v^{\ast}_k)=\la (\check{\delta} A)Y_i,v^{\ast}_k\ra$ for the standard Hopf fibration as the fibres are totally geodesic, forcing all $T$-tensor terms to vanish.
Analogous arguments apply for the Hopf fibration over $\HP^m$.
\proofend
\begin{cor}\label{cor-5.3}
The Ricci curvatures of $g_1$ are given by
\begin{align*}
&\r_1(\partial_r)=-\dim\P \frac{h''}{h}-q\frac{f''}{f};\\ &\r_1(X_i)=\dim\P\left(\frac{1-(h')^2}{h^2}\right)+\frac{2^q}{h^2}-\frac{h''}{h}-\left(\frac{h'}{h}\right)^2 -q\left(\frac{f'h'}{fh}-2\frac{f^2}{h^4}\right);\\ &\r_1(v^{\ast}_k)=(q-1)\frac{1-(f')^2}{f^2}-\frac{f''}{f}+\dim\P\left(\frac{f^2}{h^4}-\frac{f'h'}{fh}\right),
\end{align*}
with all mixed curvature terms vanishing.
\end{cor}
Notice that the metric $g_1$ extends to give a well-defined metric on $D^{n+1}$, provided $f$ and $h$ satisfy suitable boundary conditions near $r=0$.
Specifially, we require $f(0)=h(0)=0$, $f'(0)=h'(0)=1$, and $f$ and $h$ should be odd at $r=0$.
These conditions will certainly be satisfied if $f(r)=h(r)=\sin r$ for $r\in [0,\delta]$ for some small $\delta$, and we will assume this to be the case.
The values of $R$ (the radius of $D^{n+1}$) and $\delta$ will be determined later.
\begin{prp}\label{prp-5.4}
The metric $g_1=g_1(f,h)$ on $D^{n+1}$ has all Ricci curvatures strictly positive if the functions $f$ and $h$ satisfy:
\begin{align}
\label{cond-prp-5.4-1} &f(r)=h(r)=\sin r\text{ for }r\text{ small};\\
\label{cond-prp-5.4-2} &f''\le 0, h''\le 0, f''+h'' <0, f'\ge 0,\text{ and }h'\ge 0;\\
\label{cond-prp-5.4-3} &f\le h,\text{ and }\frac{f'}{f}\le\frac{h'}{h};\\
\label{cond-prp-5.4-4} &(f/h)^3\ge f'h'.
\end{align}
\end{prp}
\pf As all mixed Ricci curvature terms vanish, it suffices to show that the expressions for $\r_1(\partial_r)$, $\r_1(X_i)$ and $\r_1(v^{\ast}_k)$ are all strictly positive.
The positivity of $\r_1(\partial_r)$ is clear because of \eqref{cond-prp-5.4-2}.
To see the positivity of $\r_1(X_i)$, consider the case $q=1$.
The first term in the expression for $\r_1(X_i)$ in Corollary \ref{cor-5.3} is strictly positive for $r>0$ as a consequence of \eqref{cond-prp-5.4-1} and \eqref{cond-prp-5.4-2}.
By \eqref{cond-prp-5.4-2}, the term $-h''/h\ge 0$ for all $r$.
Therefore these two terms taken together have a strictly positive sum for all $r$.
It therefore suffices to show the non-negativity of the sum of the remaining terms:
$$2h^{-2}+(h'/h)^2-f'h'f^{-1}h^{-1}-2f^2h^{-4}.$$
By (iii) this expression is greater than or equal to
$$2h^{-2}+(h'/h)^2-(h'/h)^2-2h^{-2}$$
as required.
The case $q=3$ is analogous.
The positivity of $\r_1(v^{\ast}_k)$ follows immediately if the final term in the expression in Corollary \ref{cor-5.3} is non-negative.
But this is guaranteed by \eqref{cond-prp-5.4-4}.
\proofend
We now turn our attention to the space $G/K$.
Let $g_0$ be a bi-invariant metric on $G$ which makes $\alpha:S^q\rightarrow H/K$ an isometry, assuming the round metric of radius $1$ on $S^q$.
Note that this is possible as any bi-invariant metric on $S^1$ or $S^3$ must be round.
Fix a metric $g_{\nu}=\nu g_0$ on $G$, for some constant $\nu$.
As $g_{\nu}$ is bi-invariant it must have non-negative sectional curvature, so in particular it has non-negative Ricci curvature.
Consider the corresponding normal homogeneous metric on $G/K$.
By the O'Neill formulas, this too has non-negative sectional and therefore non-negative Ricci curvatures.
In fact, by \cite{Be} our normal homogeneous metric must have strictly positive Ricci curvature, as $\pi_1(G/K)<\infty$.
Now scale this normal homogeneous metric in the direction of the $H$-orbits by a factor $\mu$.
(Recall that $G/K$ is the total space of a fibration $H/K\hookrightarrow G/K\rightarrow G/H$.) Call the resulting metric $g_2=g_2(\mu,\nu)$.
The following is clear from the openness of the $\r>0$ condition:
\begin{obs}\label{obs-5.5}
There exists $\e_0=\e_0(G,H,K,g_0)$ such that for all $\e<\e_0$, $g_2$ has strictly positive Ricci curvature when $\mu=1+\e$.
\end{obs}
Note that $\e_0$ is independent of $\nu$. Fixing a value of $\e<\e_0$, we immediately deduce:
\begin{cor}\label{cor-5.6}
The product metric $g_1+g_2(1+\epsilon,\nu)$ on $\pro$ has positive Ricci curvature, and is both $S^q$-invariant and left $G$-invariant.
\end{cor}
The $S^q$-invariance of $g_1+g_2$ gives
\begin{cor}\label{cor-5.7}
The metric $g_1+g_2$ induces a well-defined metric $g_Q$ on $\quo$.
\end{cor}
For the purposes of gluing tubular neighbourhoods of isolated singular orbits into the manifold $M^0$, we want to ensure that the $G/K$-fibres near the boundary of $(\quo,g_Q)$ have normal homogeneous metrics.
\begin{prp}\label{prp-5.8}
Regarding $S^n$ as the total space of the Hopf fibration $S^q\hookrightarrow S^n\rightarrow\P$, consider the round metric as a submersion metric over $\P$.
Now rescale the fibres of this submersion so they are all isometric to $(S^q,\lambda ds^2_q)$.
Let $g_0$ be the bi-invariant metric on $G$ which makes $\alpha\colon S^q\rightarrow H/K$ an isometry, assuming the round metric of radius $1$ on $S^q$.
Fix a metric $g_{\nu}=\nu g_0$ on $G$, for some constant $\nu$, and consider the corresponding normal homogeneous metric on $G/K$.
Scale this normal homogeneous metric in the direction of the $H$-orbits by a factor $\mu$.
The resulting product metric on $S^n\times G/K$ induces a $G$-invariant metric on the quotient $S^n\times_{\alpha} G/K$.
The $G$-orbits in this quotient are all isometric to $G/K$ with the normal homogeneous metric induced from $g_{\nu}$ precisely when
$$\lambda=\frac{\mu\nu}{\mu-1}.$$
\end{prp}
\pf The proof is just a Cheeger-type argument, analogous to that in \cite{C}.

\proofend
\begin{cor}\label{cor-5.9}
The $G/K$-fibres at the boundary of $(\quo;g_Q)$ have normal homogeneous metrics induced by $\nu g_0$ if the function $f$ used to define $g_Q$ takes the value $(1+\e)\nu/\e$ there.
\end{cor}
We now investigate the Ricci curvature of $\quo$.
Our principal strategy for showing Ricci positivity is as follows.
If $Z_1$ and $Z_2$ are horizontal vectors in the total space of a submersion, then the Ricci curvature of the projections $\r(\check{Z}_1,\check{Z_2})$ is related to $\r(Z_1,Z_2)$ by the following O'Neill formula \cite{B}, \S9.36c:
$$\r(Z_1,Z_2)=\r(\check{Z}_1,\check{Z_2})-2\la A_{Z_1},A_{Z_2}\ra -\la TZ_1,TZ_2\ra +\frac{1}{2}\left(\la\nabla_{Z_1}N,Z_2\ra +\la\nabla_{Z_2}N,Z_1\ra\right).$$
In particular, this means that for $Z=Z_1=Z_2$ we have
$$\r(\check{Z})=\r(Z)+2\la A_Z,A_Z\ra +\la TZ,TZ\ra -\la\nabla_Z N,Z\ra.$$
Clearly $\la A_Z,A_Z\ra\ge 0$ and $\la TZ,TZ\ra\ge0$.
Therefore, assuming $\r(Z)>0$, if $\la\nabla_Z N,Z\ra\le 0$ we must have $\r(\check{Z})>0$ also.
It is therefore crucial to understand the $\la\nabla_Z N,Z\ra$ term.
In order to do this, we must first identify the vector field $N$.
Recall that by definition (\cite{B} \S9.34), $N=\sum_k T_{U_k}U_k$, where $\{U_k\}$ is an orthonormal frame field for the vertical distribution.
Although we will not need to know such a vertical frame field explicitly for the computation of $N$, we will need such formulas later on.
With this in mind, set
$$U_k=\frac{1}{{\sqrt{\den}}}(v_k^{\ast},-v^{\alpha}).$$
\begin{lm}\label{lm-5.10}
For the metric $g_1+g_2$ on $\pro$ we have
$$N=\frac{-qff'}{\den}\partial_r.$$
\end{lm}\pf We can view the metric $g_1+g_2$ on $\pro$ as a submersion metric which has been created from a submersion with isometric, totally geodesic fibres by rescaling in fibre directions by $\den$.
All $T$-tensor terms for the totally geodesic submersion vanish.
The effect of the rescaling on all the quantities appearing the O'Neill formulas was computed in \cite{W1}.
The expression for $N$ following a metric rescale by a function $\theta$ defined on the base is given by
$$N=-\frac{\dim(\text{fibre})}{2\theta}\nabla\theta.$$
Setting $\theta=\den$ in this formula gives the desired expression.
\proofend
For convenience we will write $N=\phi\partial_r$ from now on, with
$$\phi=\phi(r)=\frac{-qff'}{\den}.$$ \par
In the following, metric quantities without subscript will refer to the product metric $g_1+g_2$ on $\pro$.
Recall that $g_Q$ is the induced metric on the quotient space $\quo$.
It is easy to see that the horizontal distribution in $\pro$ is
$$(\mathcal{H}\oplus 0)\oplus\left\{\left(\frac{\mu\nu}{\lambda} v^{\ast},v^{\alpha}\right) \,|\, v\in TS^q\right\}\oplus (0\oplus \mathfrak{m}),$$
where $\mathcal{H}$ denotes the horizontal distribution for the Hopf submersion metric on $S^n$, and where $\mathfrak{m}$ is the distribution of orthogonal complements to $H/K$ orbits in $G/K$.
Recall that the vector fields $\{X_i\}$ are an orthonormal basis for $\mathcal{H}\oplus 0$.
Let $\{w_j\}$ be an orthonomal frame field for $\mathfrak{m}$.
Setting $\Delta_k=\vk$, we have that $\{\Delta_k\}$ is an orthogonal frame field for $\{((\mu\nu/\lambda) v^{\ast},v^{\alpha})\}$.
However, note that $\{\Delta_k\}$ is not an orthonormal set as
$$\|\Delta_k\|^2=\frac{(1+\e)\nu(f^2+(1+\e)\nu)}{f^2}.$$
As before, the projections of any of these vectors to $\quo$ will be indicated by a \ $\check{}$\ .
Collectively, these projections form a local basis. \par
Although the vectors $\{X_i\}$ are orthonormal, it will sometimes be useful in subsequent calculations to write $X_i=hY_i$, so the projections of the $Y_i$ on the base are unit vector fields with respect to the Fubini-Study metric.
The significance of this is that $\{Y_i\}$ are independent of the $r$ parameter.
As a result,
$$[Y_i,\partial_r]=[Y_i,\Delta_k]=[Y_i,w_j]=[Y_i,N]=[Y_i,U_k]=0.$$
\begin{lm}\label{lm-5.11} The following formulas hold:
\begin{align*}
[X_i,N]&=\phi(h'/h)X_i \ ; \\
[\partial_r,N]&=\phi' \partial_r \ ; \\
[v^{\ast}_k,N]&=[v^{\alpha},N]=[w_j,N]=0 .
\end{align*}
\end{lm}
\pf For the first expression, we begin by writing the Lie bracket as $[hY_i,\phi\partial_r]$ and then using the fact that $[Y_i,\partial_r]=0$.
The second expression is an elementary calculation, and the vanishing of the final three terms is immediate since $v^{\ast}$, $v^{\alpha}$ and $w_j$ are all independent of $r$ and are tangent to different factors to $(0,R]$ in the product $(0,R]\times S^n\times G/K$.
\proofend
\begin{lm}\label{lm-5.12}
The vector field $N$ has the following covariant derivatives:
$$\nabla_{X_i} N=\phi\frac{h'}{h} X_i, \hskip 2cm\nabla_{w_j}N=0, \hskip 2cm\nabla_{\partial_r}N=\phi'\partial_r,$$
$$\nabla_{v_k^{\ast}}N=\phi f'f^{-1}v_k^{\ast}, \quad\nabla_{\Delta_k}N=\phi f'f^{-1}\frac{(1+\e)\nu}{(1+\e)\nu+f^2}\Delta_k.$$
\end{lm}
\pf We proceed using the Koszul formula:
$$2\la\nabla_A B,C\ra=A\la B,C\ra+B\la C,A\ra-C\la A,B\ra +[[A,B],C]-[[B,C],A]+[[C,A],B].$$
We compute each of the sixteen possible terms $\la\nabla_{\bullet} N,\star\ra$.
For the most part these are zero.
We briefly mention those which are not. \par
For $\nabla_{X_i} N$, the term $2\la\nabla_{X_i} N,X_j\ra =\la [X_i,N],X_j\ra -\la [N,X_j],X_i\ra$.
By Lemma \ref{lm-5.11} we see that each of the terms on the right-hand side is equal to $\phi(h'/h)X_i$, and hence $\la\nabla_{X_i} N,X_j\ra = \phi\frac{h'}{h}\delta_{ij}$. \par
For $\nabla_{\partial_r}N$, the term $2\la\nabla_{\partial_r}N, \partial_r\ra =\la [\partial_r,N],\partial_r\ra -\la [N, \partial_r], \partial_r\ra$.
By Lemma \ref{lm-5.11} both of these terms are equal to $\phi'\partial_r$, and hence $\la\nabla_{\partial_r}N, \partial_r\ra =\phi'$. \par
For $\nabla_{v_k^{\ast}}N$ the only non-zero expression is $2\la\nabla_{v_k^{\ast}}N, v_l^{\ast}\ra=N\la v_k^{\ast},v_l^{\ast}\ra$.
Now $\|v_k^{\ast}\|^2=f^2$, so $\la v_k^{\ast},v_l^{\ast}\ra=f^2\delta_{kl}$, giving $\la\nabla_{v_k^{\ast}}N, v_l^{\ast}\ra=\phi ff' \delta_{kl}$. \par
For $\nabla_{\Delta_k} N$, the only non-zero term is $\la\nabla_{\Delta_k} N,\Delta_k\ra$.
We have
\begin{align*}
2\la\nabla_{\Delta_k} N,\Delta_k\ra &= 2\la\nabla_{(1+\e)\nu f^{-2}v_k^{\ast}+v_k^{\alpha}}N,(1+\e)\nu f^{-2}v_l^{\ast}+v_l^{\alpha}\ra\\
&=2\left[\fra\right]^2\la\nabla_{v_k^{\ast}}N, v_l^{\ast}\ra \text{ as all other terms clearly vanish},\\ &=2\left[\fra\right]^2 \phi ff' \delta_{kl} \text{ as shown above}.
\end{align*}
Therefore $\la\nabla_{\Delta_k} N,\Delta_k\ra=(1+\e)^2\nu^2\phi f'/f^3$. \par
From the above results, the conclusion of the Lemma is easy to establish.
Note that for $\nabla_{\Delta_k}N$ we need to take care since $\Delta_k$ is not a unit vector.
Specifically, we have
\begin{align*}
\nabla_{\Delta_k}N&=(1+\e)^2\nu^2\phi\frac{f'}{f^3}\frac{1}{\|\Delta_k\|}\frac{\Delta_k}{\|\Delta_k\|}\\
&=\phi\frac{f'}{f}\frac{(1+\e)\nu}{(1+\e)\nu+f^2}\Delta_k.
\end{align*}
\proofend
\begin{cor}\label{cor-5.13}
For the metric $g_Q$ on $\quo$, we have $\r_Q(\check{w}_j)>0$ for each $j$.
\end{cor}
\pf By Observation \ref{obs-5.5} we have $\r(w_j)>0$ for all $j$.
Now $\r(w_j)$ and $\r_Q(\check{w}_j)$ are related by
$$\r_Q(\check{w}_j)=\r(w_j)+2\la A_{w_j},A_{w_j}\ra+\la Tw_j,Tw_j\ra -\la\nabla_{w_j}N,w_j\ra.$$
Thus the first three terms on the right-hand side are strictly positive, non-negative and non-negative respectively.
By Lemma \ref{lm-5.12} we see that the final term vanishes, which establishes the result.
\proofend
\begin{cor}\label{cor-5.14}
For every $i$ we have $\r_Q(\check{X}_i)>0$.
\end{cor}
\pf We argue as in the proof of Corollary \ref{cor-5.13} above.
The only difference this time is that the term $\la\nabla_{X_i} N,X_i\ra$ is not zero.
From Lemma \ref{lm-5.12} we have $\la\nabla_{X_i} N,X_i\ra = \phi\frac{h'}{h}$.
But recall that
$$\phi=-\frac{qff'}{\den}\partial_r,$$
so in particular we have $\phi(r)\le 0$ for all $r$.
This means that the $\la\nabla_{X_i} N,X_i\ra$ term is non-positive, and therefore makes a non-negative contribution to $\r_Q(\check{X}_i)$.
\proofend
\begin{cor}\label{cor-5.15}
For every $k$ we have $\r_Q(\check{\Delta}_k)>0$.
\end{cor}
\pf The proof is essentially the same as for Corollary \ref{cor-5.14}.
This time we have
$$\la\nabla_{\Delta_k} N,\Delta_k\ra=(1+\e)^2\nu^2\phi\frac{f'}{f^3}$$
from Lemma \ref{lm-5.12}, and this makes a non-negative contribution to $\r_Q(\check{\Delta}_k)>0$.
\proofend
\begin{lm}\label{lm-5.16}
We have $\la T\partial_r,T\partial_r\ra= q(f'f)^2[\den]^{-2}$.
All other terms of the form $\la T\bullet,T\star\ra$ vanish.
\end{lm}
\pf We can view the metric $g_1+g_2$ on $\pro$ as a submersion metric which has been created from a submersion with isometric, totally geodesic fibres by rescaling in fibre directions by $\den$.
All $T$-tensor terms for the totally geodesic submersion vanish.
The effect of this rescaling on all the quantities appearing the O'Neill formulas was computed in \cite{W1}.
For the term $\la TA,TB\ra$, the value following a metric rescale by a function $\theta$ defined on the base is given by
$$\frac{\dim(\text{fibre})}{4\theta^2}A(\theta)B(\theta).$$
Setting $\theta=\den$ and computing derivatives gives the result.
\proofend
\begin{lm}\label{lm-5.17}
For each $k$ we have
$$A_{\partial_r} U_k=\frac{ff'}{[\den]^{\frac{3}{2}}}\Delta_k.$$
\end{lm}
\pf As the $A$-tensor is linear in both entries,
$$A_{\partial_r} U_k=\frac{1}{\sqrt{\den}}A_{\partial_r} (v_k^{\ast},-v^{\alpha}).$$
By definition of the $A$-tensor, this quantity is
$$\frac{1}{\sqrt{\den}}\mathcal{H}\nabla_{\partial_r}(v_k^{\ast},-v^{\alpha}),$$
where $\mathcal{H}$ denotes the horizontal component. \par
We next compute the components of this vector in the various directions.
Using the formula for $\nabla_{v_k^{\ast}}N$ established in Lemma \ref{lm-5.12} we see that
\begin{align*}
\la\nabla_{\partial_r}(v_k^{\ast},-v^{\alpha}),\Delta_l\ra &=f'f^{-1}\la v_k^{\ast}, \Delta_l\ra\\ &=(1+\e)\nu f'f^{-1}\delta_{kl}.
\end{align*}
Similarly, for the other directions it is easy to see that the terms $\la\nabla_{\partial_r}(v_k^{\ast},-v^{\alpha}),X_i\ra$, $\la\nabla_{\partial_r}(v_k^{\ast},-v^{\alpha}),w_j\ra$ and $\la\nabla_{\partial_r}(v_k^{\ast},-v^{\alpha}),\partial_r\ra$ all vanish. \par
Bearing in mind the fact that $\Delta_k$ is not a unit vector, we deduce that
\begin{align*}
\nabla_{\partial_r}(v_k^{\ast},-v^{\alpha})&=(1+\e)\nu\frac{f'}{f}\frac{1}{\|\Delta_k\|}\frac{\Delta_k}{\|\Delta_k\|}\\
&=\frac{ff'}{[\den]^\frac{3}{2}}\Delta_k
\end{align*}
as required.
\proofend
\begin{cor}\label{cor-5.18}
$$\la A_{\partial_r},A_{\partial_r}\ra= q\frac{(1+\e)\nu (f')^2}{[\den]^2}.$$
\end{cor}
\pf Recall from \cite{B}, \S9.33 that
$$\la A_{\partial_r},A_{\partial_r}\ra=\sum_k\la A_{\partial_r}U_k,A_{\partial_r}U_k\ra.$$
Using the result of Lemma \ref{lm-5.17}, we immediately obtain the desired expression.
\proofend
We are now in a position to investigate $\r_Q (\check{\partial}_r)$.
\begin{prp}\label{prp-5.19}
$\r_Q (\check{\partial}_r)>0$.
\end{prp}
\pf We know that
$$\r_Q (\check{\partial}_r)=\r(\partial_r)+2\la A_{\partial_r},A_{\partial_r}\ra +\la T\partial_r,T\partial_r\ra -\la\nabla_{\partial_r}N,\partial_r\ra.$$
Using the formula for $\r(\partial_r)$ from Corollary \ref{cor-5.3}, the formulas of Corollary \ref{cor-5.18} and Lemma \ref{lm-5.16} for the next two terms on the right-hand side, and Lemma \ref{lm-5.12} to evaluate the final term we obtain the expression
\begin{align*}
\r_Q(\check{\partial}_r)=-&\dim\P\frac{h''}{h}-q\frac{f''}{f}+q\frac{(1+\e)\nu (f')^2}{[\den]^2}+ q\frac{(ff')^2}{[\den]^2}\cr\cr &+q\frac{(f')^2}{\den}+q\frac{ff''}{\den}-2q\frac{(ff')^2}{[\den]^2},
\end{align*}
where the final three terms are just $\phi'$ written out explicitly. \par
Collecting similar terms gives
$$-\dim\P\frac{h''}{h}-q\frac{f''}{f}\left[1-\frac{f^2}{\den}\right]+ q\frac{(ff')^2}{[\den]^2}\left[1+\frac{2(1+\e)\nu}{f^2}+\frac{\den}{f^2}-2\right].$$
Simplifying this gives
$$r_Q(\partial_r)=-\dim\P\frac{h''}{h}+q\frac{(1+\e)\nu}{\den}\left(\frac{3(f')^2}{\den}-\frac{f''}{f}\right).$$
As we are assuming that both $f$ and $h$ are concave down functions, at least one of which is strictly concave down for all $r$, we deduce that this expression is strictly positive as claimed.
\proofend
So far we have established that $\r_Q(\check{\partial}_r)$, $\r_Q(\check{X}_i)$, $\r_Q(\check{\Delta}_k)$ and $\r_Q(\check{w}_j)$ are all strictly positive.
However this is not sufficient to deduce that {\it all} Ricci curvatures of the metric $g_Q$ are strictly positive.
\begin{lm}\label{lm-5.20}
For any $a,b,c\in\R$, $\r_Q(a\check{X}_i+b\check{\Delta}_k+c\check{w}_j)>0$.
\end{lm}
\pf Using the same line of reasoning as employed in Corollaries \ref{cor-5.13}, \ref{cor-5.14} and \ref{cor-5.15}, it suffices to show that the expression $\la\nabla_{aX_i+b\Delta_k +c\check{w}_j} N,aX_i+b\Delta_k+c\check{w}_j\ra$, is non-positive.
But this follows easily from (the proof of) Lemma \ref{lm-5.12}.
\proofend
It remains to study Ricci curvatures of the form $\r_Q(\check{\partial}_r+\check{Z})$ for $Z\in\Span\{X_i\}\oplus\Span\{\Delta_k\}\oplus\Span\{w_j\}$.
By elementary linear algebra, Ricci curvatures of this form will be positive if and only if
$$\r_Q(\check{\partial}_r)\r_Q(\check{Z})>(\r_Q(\check{\partial}_r,\check{Z}))^2$$
for all $Z$.
\begin{prp}\label{prp-5.21}
For all $Z\in\Span\{X_i\}\oplus\Span\{\Delta_k\}\oplus\Span\{w_j\}$, we have
$$\r_Q(\check{\partial}_r)\r_Q(\check{Z})>(\r_Q(\check{\partial}_r,\check{Z}))^2.$$
\end{prp}
\pf We prove this proposition in several steps.
The first step is to establish the inequality
\[{pf-prp-5.21-d1}\r_Q(\check{\partial}_r)\ge 2\la A_{\partial_r},A_{\partial_r}\ra.\]
From the curvature formulas established in the proof of Proposition \ref{prp-5.19}, we see that our inequality is equivalent to
$$-\dim\P\frac{h''}{h}+q\frac{(1+\e)\nu}{\den}\left(\frac{3(f')^2}{\den}-\frac{f''}{f}\right)\ge 2q\frac{(1+\e)\nu(f')^2}{[\den]^2}.$$
As $h''/h\ge 0$ it suffices to show that
$$\frac{3(f')^2}{\den}-\frac{f''}{f}\ge\frac{2(f')^2}{\den}.$$
As $f''/f\ge 0$ it then suffices to show that
$$\frac{3(f')^2}{\den}\ge\frac{2(f')^2}{\den},$$
which is clearly true.
Thus \eqref{pf-prp-5.21-d1} is established. \par
In fact we can go further than this.
By assumption, at least one of $f$ or $h$ is strictly concave down for each $r$.
Thus the inequality \eqref{pf-prp-5.21-d1} can actually be replaced by
\[{pf-prp-5.21-d2}\r_Q(\check{\partial}_r) > 2\la A_{\partial_r},A_{\partial_r}\ra.\] \par
We next claim that for any $Z\in\Span\{X_i\}\oplus\Span\{\Delta_k\}\oplus\Span\{w_j\}$, we have
$$\r_Q(\check{Z})> 2\la A_Z,A_Z\ra.$$
To see this, note that by the O'Neill formulas,
$$\r_Q(\check{Z})-\r(Z)-\la TZ,TZ\ra +\la\nabla_Z N,Z\ra =2\la A_Z,A_Z\ra.$$
Now we know from Corollary \ref{cor-5.6} that $\r(Z)>0$, $\la TZ,TZ\ra=0$ by Lemma \ref{lm-5.16}, and $\la\nabla_Z N,Z\ra\le 0$ by Lemma \ref{lm-5.12}.
Thus the inequality follows. \par
It is an elementary consequence of the Cauchy-Schwarz inequality that
$$\la A_Z,A_{\partial_r}\ra^2\le\la A_Z,A_Z\ra\la A_{\partial_r},A_{\partial_r}\ra,$$
and in particular we have
$$4\la A_Z,A_{\partial_r}\ra^2\le [2\la A_Z,A_Z\ra] [2\la A_{\partial_r},A_{\partial_r}\ra].$$
Combining this with \eqref{pf-prp-5.21-d2} and the corresponding inequality for $\r_Q(\check{Z})$ gives
$$\r_Q(\check{Z}) \r_Q(\check{\partial}_r)>4\la A_Z,A_{\partial_r}\ra^2.$$ \par
The proof of the Proposition will now follow from our final claim: $\r_Q(\check{Z},\check{\partial}_r)=2\la A_Z,A_{\partial_r}\ra$.
To see this we use the O'Neill formula
$$\r_Q(\check{Z},\check{\partial}_r)=\r(Z,\partial_r)+ 2\la A_Z,A_{\partial_r}\ra+\la TZ,T\partial_r\ra -\frac{1}{2}\left[\la\nabla_Z N,\partial_r\ra +\la\nabla_{\partial_r} N,Z\ra\right].$$
By Lemma \ref{lm-5.16} we have $\la TZ,T\partial_r\ra=0$, and the vanishing of the final term follows from Lemma \ref{lm-5.12}.
Thus the claim, and hence the Proposition is established.
\proofend
We immediately deduce:
\begin{cor}\label{cor-5.22}
 For all $Z\in\Span\{X_i\}\oplus\Span\{\Delta_k\}\oplus\Span\{w_j\}$, we have $\r_Q(\check{\partial}_r+\check{Z})>0$.
\end{cor}
We are now in a position to prove Theorem \ref{tube}. \par
\pfo {Theorem \ref{tube}.}
We show that the functions $f$ and $h$ can be chosen so that the metric $g_Q$ satisfies all the requiements of the Theorem.
By combining the results of Corollaries \ref{cor-5.13}, \ref{cor-5.14} and \ref{cor-5.15}, Proposition \ref{prp-5.19}, Lemma \ref{lm-5.20} and Corollary \ref{cor-5.22}, we see that all Ricci curvatures of the metric $g_Q$ are strictly positive, provided $f$ and $h$ satisfy the conditions of Proposition \ref{prp-5.4}.
Of course $g_Q$ is $G$-invariant by construction.
By Corollary \ref{cor-5.9}, we will obtain fibres in an $\iota$-neighbourhood of the boundary all isometric to the normal homogeneous metric induced by $\nu g_0$ (for any choice of $\nu$) if $f(r)=(1+\e)\nu/\e$ for all $r\in [R-\iota,R]$.
Choose a function $h(r)$ such that $h(r)=\sin r$ for $r$ small, $h''(r)<0$ for all $r$, and $h(R)=\lambda$, $h'(R)=\Lambda$ for some $R>0$.
It is clear that we can make such a choice.
Next, note that if we set $f(r)=h(r)$, the conditions laid out in Proposition \ref{prp-5.4} are all satisfied.
However, setting $f(r)=h(r)$ will not allow us to achieve the required $f'(r)=0$ for $r\in [R-\iota,R]$.
For any choice of $\delta$ such that $0<\delta<R$, let $f_0(r)=h(r)$ for $r\in [0,\delta]$ and $f_0(r)=h(\delta)$ for $r\in (\delta,R]$.
The function $f_0$ is clearly not smooth, however it is clear that we can smooth it in an arbitrarily small neighbourhood of $r=\delta$ to a function $f$, so that $f''\le 0$.
Provided the smoothing neighbourhood is sufficiently small, the functions $f$ and $h$ then satisfy the requirements of Proposition \ref{prp-5.4}.
Note that we can arrange for $f(R)$ to be any value less than $h(R)=\lambda$.
To complete the proof, it remains to show that $f$ can take the value $(1+\e)\nu/\e$ close to the boundary.
But $f$ can take any value less than $\lambda$ at the boundary, so provided $(1+\e)\nu/\e<\lambda$ this boundary condition can be achieved.
Rearranging, this gives $\nu<\lambda\e/(1+\e)$ as claimed.
\proofend

\section{Open problems}\label{sec-open}
We conclude the paper with a selection of geometric open problems.
\begin{ils}
\item\label{op-1} Do any manifolds with a single singular orbit and cohomogeneity greater than three admit an invariant metric with positive Ricci curvature?
Recall that by Theorem \ref{one} we can construct manifolds of cohomogeneity three with a single singular orbit and positive Ricci curvature.
The problem with extending this family into higher cohomogeneities is that it necessitates extending the Fubini-Study metric on $\CP^{(k-1)/2}$ over the disc bundle corresponding to the imaginary sub-bundle of the canonical quaternionic line bundle over $\HP^{(k-3)/4}$ in such a way that the extension satisfies the requirements of Theorem \ref{thm-riccurvmain}. It is not clear to the authors whether such an extension is possible.

\item If the answer to question \ref{op-1} is yes, then do any of these manifolds admit invariant metrics with non-negative sectional curvature?
\item The Ricci positive examples displayed in section 2 have at most two singular orbits.
Is it possible to find invariant Ricci positive metrics on manifold having more than two singular orbits? \smallskip \noindent The obvious candidates are those for which $B$ (the space of orbits obtained when tubular neighbourhoods of the singular orbits have been removed from the original manifold) is a $3$-sphere less some discs.
(Thus the boundary is a disjoint union of $2$-spheres, that is, $\CP^1$s.) It is easily checked that conditions \eqref{cond-fbs} and \eqref{cond-pricurv} of Theorem \ref{thm-riccurvmain} mean that while two discs can comfortably be removed, taking out three discs results in these conditions {\it just} failing to hold.
The same is true when $B$ is a $5$-sphere less some discs (so the boundary components are all equal to $S^4=\HP^1$.) It is not clear whether the failure of these obvious candidates is due to their special nature, or whether they represent a general phenomenon.
Indeed it might be possible that {\it no} manifold with more than two singular orbits can support an invariant metric with positive Ricci curvature.
\item Are there any simply-connected examples which do {\it not} admit an invariant metric with positive Ricci curvature? \smallskip \noindent As noted in the Introduction, there is very little chance of all simply-connected $G$-manifolds with finitely many singular orbits admitting invariant metrics with positive Ricci curvature, since the topology of the space of orbits can be highly non-trivial and must surely influence the possible curvatures which the manifold can display.
\item\label{op-5} Are there any examples of non-double manifolds among the families $M^{11}_{p_1p_2q_1q_2}$ or $M^{13}_{p_1p_2q_1q_2}$, or indeed any non-double examples of any kind, which admit invariant metrics of non-negative sectional curvature? \smallskip \noindent If we simply want to join tubular neighbourhoods of two different singular orbits to create our manifold, then the main problem is that the horizontal distributions arising from the metric construction process never seem to match.
A possible strategy here is to look for horizontal distributions which could be deformed so as to join smoothly, whilst preserving non-negative sectional curvature.
The authors have no idea when or how such a deformation might be possible.
\item Do any of the manifolds with two different singular orbits admit metrics of almost non-negative sectional curvature? \smallskip \noindent On the face of it, this question is more likely to have a positive answer than question \ref{op-5}.
The motivation for this question arises from \cite{ST}, where it is shown that every compact cohomogeneity one manifold admits such a metric.
One of the features of cohomogeneity one manifolds which is important here is that the space of orbits, being one-dimensional, makes no contribution to the curvature.
On the other hand, in our situation this is not the case, so such metrics will almost certainly be much less common (assuming they exist at all).
\end{ils}

\bigskip
\noindent Department of Mathematics and Statistics, National University of Ireland Maynooth, Maynooth, Co. Kildare, Ireland.
Email: stefan@maths.nuim.ie, david.wraith@nuim.ie.

\begin{thebibliography}{10}
\bibitem{B} A.L. Besse, {\it Einstein Manifolds}, Springer-Verlag, Berlin (2002).
\bibitem{Be} V. Berestovskii, {\it Homogeneous Riemannian manifolds of positive Ricci curvature}, Math Notes {\bf 58} no. 3 (1995), 905-909.
\bibitem{BW0} S. Bechtluft-Sachs, D. J. Wraith, {\it Manifolds of low cohomogeneity and positive Ricci curvature}, Diff. Geom. Appl. {\bf 28} (2010), 282-289.
\bibitem{BW} S. Bechtluft-Sachs, D. J. Wraith, {\it On the topology of $G$-manifolds with finitely many non-principal orbits}, arXiv:1106.3432.
\bibitem{Br} G. Bredon, {\it Introduction to compact transformation groups}, Pure and Applied Mathematics {\bf 46}, Academic Press, New York-London (1972).
\bibitem{C} J. Cheeger, {\it Some examples of manifolds of nonnegative curvature}, J. Differential Geometry, {\bf 8} (1973), 623-628.
\bibitem{D} O. Dearricott, {\it A 7-manifold with positive curvature}, preprint (2008).
\bibitem{GVWZ} K. Grove, L. Verdiani, B. Wilking, W. Ziller, {\it Non-negative curvature obstructions in cohomogeneity one and the Kervaire spheres}, Ann. Sc. Norm. Super. Pisa Cl. Sci. (5) {\bf 5} (2006) no. 2, 159-170.
\bibitem{GVZ} K. Grove, L. Verdiani, W. Ziller {\it A new type of a positively curved manifold}, arXiv:0809.2304v2 [math.DG].
\bibitem{GWZ} K. Grove, B. Wilking, W. Ziller, {\it Positively curved cohomogeneity-one manifolds and 3-Sasakian geometry}, J. Diff. Geom. {\bf 78} (2008) no. 1, 33-111.
\bibitem{GZ1} K. Grove, W. Ziller, {\it Curvature and symmetry of Milnor spheres}, Ann. of Math. (2) {\bf 152} (2000), no. 1, 331--367.
\bibitem{GZ2} K. Grove, W. Ziller, {\it Cohomogeneity one manifolds with positive Ricci curvature}, Invent. Math. {\bf 149} (2002), 619-646.
\bibitem{H} B. Hanke, {\it Positive scalar curvature with symmetry}, arXiv:math.GT/0512284.
\bibitem{LY} H. B. Lawson, S.-T. Yau, {Scalar curvature, non-abelian group actions and the degree of symmetry of exotic spheres}, Comment. Math. Helv. {\bf 49} (1974), 232-244.
\bibitem{P} G. Perelman, {\it Construction of manifolds of positive Ricci curvature with big volume growth and large Betti numbers}, $\lq$Comparison Geometry', Cambridge University Press, (1997).
\bibitem{PW} P. Petersen, F. Wilhelm, {\it An exotic sphere with positive sectional curvature}, arXiv:0805.0812.
\bibitem{ST} L. Schwachh\"ofer, W. Tuschmann, {\it Almost non-negative curvature and cohomogeneity one}, preprint no. 62/2001, Max-Planck-Institut Leipzig, (2001).
\bibitem{V1} L. Verdiani, {\it Cohomogeneity one Riemannian manifolds of even dimension with strictly positive sectional curvature}, Math. Z. {\bf 241} (2002), 329-339.
\bibitem{V2} L. Verdiani, {\it Cohomogeneity one manifolds of even dimension with strictly positive sectional curvature}, J. Diff. Geom. {\bf 68} (2004), 31-72.
\bibitem{V} J. Vilms, {\it Totally geodesic maps}, J. Differential Geometry {\bf 4} (1970), 73-79.
\bibitem{W1} D. J. Wraith, {\it Exotic spheres with positive Ricci curvature}, PhD Thesis, University of Notre Dame, 1996.
\bibitem{W2} D. J. Wraith, {\it Bundle stabilisation and positive Ricci curvature}, Diff. Geom. Appl. {\bf 25} (2007), 552-560.
\end{thebibliography}
\end{document}